\documentclass[11pt]{amsart}

\usepackage[dvips]{epsfig}
\usepackage{graphics}
\usepackage{latexsym}
\usepackage{verbatim}
\usepackage{amsmath}
\usepackage{amsthm}
\usepackage{amssymb}
\newtheorem{theorem}{Theorem}[section]
\newtheorem{lemma}[theorem]{Lemma}

\def\Remark{\medskip\noindent{\bf Remark: }}
\def\Remarks{\medskip\noindent{\bf Remarks: }}

\newcommand{\ens}[1]{\mathbb{#1}}

\newcommand{\N}{\mathbb{N}}

\newcommand{\R}{\mathbb{R}}

\def\cal{\mathcal}

\def\derpar#1#2{\frac{\partial#1}{\partial#2}}
\def\var{\varepsilon}
\def\signcm{\bigskip\bigskip\hspace{80mm}
\vbox{{\sc C. Mouhot\par\vspace{3mm}
UMPA, ENS Lyon\par
46 all\'ee d'Italie\par
69364 Lyon Cedex 07\par
FRANCE\par\vspace{3mm}
e-mail:} cmouhot@umpa.ens-lyon.fr }}

\def\signlp{\bigskip\bigskip\hspace{80mm}
\vbox{{\sc L. Pareschi\par\vspace{3mm} Universit\`{a} di
Ferrara\par Via Machiavelli 35\par I-44100 Ferrara \par ITALY
\par\vspace{3mm} e-mail:} pareschi@dm.unife.it }}

\begin{document}

\title[Fast algorithms for computing the Boltzmann collision operator]
{Fast algorithms for computing the Boltzmann collision operator}

\author{Cl\'ement Mouhot and Lorenzo Pareschi}

\hyphenation{bounda-ry rea-so-na-ble be-ha-vior pro-per-ties
cha-rac-te-ris-tic}

\begin{abstract}
The development of accurate and fast numerical schemes for the
five fold Boltzmann collision integral represents a challenging
problem in scientific computing. For a particular class of
interactions, including the so-called {\em hard spheres model} in
dimension three, we are able to derive spectral methods
that can be evaluated through fast
algorithms. These algorithms are based on a suitable
representation and approximation of the collision operator.
Explicit expressions for the errors in the schemes are given
and spectral accuracy is
proved. Parallelization properties and adaptivity of the
algorithms are also discussed.
\end{abstract}

\maketitle



\tableofcontents

\section{Introduction}\label{sec:mp:intro}
\setcounter{equation}{0}



The Boltzmann equation describes the behavior of a dilute gas of
particles when the only interactions taken into account are binary
elastic collisions. It reads for $x,v \in \R^d$ ($d \ge 2$)
  \begin{equation*}
  \derpar{f}{t} + v \cdot \nabla_x f = Q(f,f)
  \end{equation*}
where $f(t,x,v)$ is the time-dependent particle distribution
function in the phase space. The Boltzmann collision operator $Q$
is a quadratic operator local in $(t,x)$. The time and position
acts only as parameters in $Q$ and therefore will be omitted in
its description
  \begin{equation}\label{eq:mp:Q}
  Q (f,f)(v) = \int_{\R^d \times\ens{S}^{d-1}} B(|v-v_*|,\cos \theta) \,
  \left( f'_* f' - f_* f \right) \, dv_* \, d\sigma.
  \end{equation}
In~\eqref{eq:mp:Q} we used the shorthand $f = f(v)$, $f_* = f(v_*)$,
$f ^{'} = f(v')$, $f_* ^{'} = f(v_* ^{'})$. The velocities of the
colliding pairs $(v,v_*)$ and $(v',v'_*)$ are related by
  \begin{equation*}
  v' = \frac{v+v_*}{2} + \frac{|v-v_*|}{2} \sigma, \qquad
  v'_* = \frac{v+v^*}{2} - \frac{|v-v_*|}{2} \sigma.
  \end{equation*}
The collision kernel $B$ is a non-negative function which by
physical arguments of invariance only depends on $|v-v_*|$ and
$\cos \theta = {\hat g} \cdot \sigma$ (where ${\hat g} =
(v-v_*)/|v-v_*|$).

Boltzmann's collision operator has the fundamental properties of
conserving mass, momentum and energy
  \begin{equation*}
  \int_{\R^d} Q(f,f)\phi(v) \,dv = 0, \quad
  \phi(v)=1,v_1, \dots, v_d,|v|^2
 \end{equation*}
and satisfies the well-known Boltzmann's $H$ theorem
  \begin{equation*} 
  - \frac{d}{dt} \int_{\R^d} f \log f \, dv = - \int_{\R^d} Q(f,f) \log(f) \, dv \geq 0.
  \end{equation*}
The functional $- \int f \log f$ is the entropy of the solution.
Boltzmann's $H$ theorem implies that any equilibrium
distribution function, {\em i.e.} any function which is a maximum of the entropy,
has the form of a locally Maxwellian distribution
  \begin{equation*}
  M(\rho,u,T)(v)=\frac{\rho}{(2\pi T)^{d/2}}
  \exp \left\{ - \frac{\vert u - v \vert^2} {2T} \right\},
  \end{equation*}
where $\rho,\,u,\,T$ are the density, mean velocity
and temperature of the gas, defined by
  \begin{equation*}
  \rho = \int_{\R^d}f(v) \, dv, \quad
  u = \frac{1}{\rho}\int_{\R^d} v f(v) \, dv, \quad
  T = {1\over{d\rho}} \int_{\R^d}\vert u - v \vert^2f(v) \, dv.
  \end{equation*}
For further details on the physical background and derivation of
the Boltzmann equation we refer to~\cite{CIP:94,Vill:hand}.

The construction of numerical methods for Boltzmann equations
represents a real challenge for scientific computing and it is of
paramount importance in many applications, ranging from {rarefied
gas dynamics} (RGD)~\cite{CIP:94}, plasma physics~\cite{DeLu:92},
granular flows~\cite{BCP:97,BCP:erra:97}, semiconductors~\cite{MRS:89} and
quantum kinetic theory~\cite{EsMi:01}.

Most of the difficulties are due to the multidimensional structure
of the collisional integral $Q$, as the integration runs
on a $5$-dimensional unflat manifold.
In addition to the unpracticable
computational cost of deterministic quadrature rules the
integration has to be handled carefully since it is at the basis
of the macroscopic properties of the equation. Additional
difficulties are represented by the stiffness induced by the
presence of small scales, like the case of small mean free
path~\cite{GaPaTo:97} or the case of large
velocities~\cite{FiPa:02}.

For such reasons realistic numerical computations are based on
probabilistic Monte-Carlo techniques at different levels. The most
famous examples are the {direct simulation Monte Carlo (DSMC)}
methods by Bird~\cite{Bird:94} and by Nanbu~\cite{Nanb:83}. These methods
preserve the conservation properties of the equation in a natural
way and avoid the computational complexity of a deterministic
approach. However avoiding the low accuracy and the fluctuations
of the results becomes extremely expensive in presence of
nonstationary flows or close to continuum regimes.

Among deterministic approximations, one of the most popular
methods in RGD is represented by the {discrete velocity models}
(DVM) of the Boltzmann equation. These methods~\cite{Buet:96,
MaSc:FBE:92, BoVaPaSc:cons:95, CoRoSc:hom:92, HePa:DVM:02,
RoSc:quad:94} are based on a regular grid in the velocity field
and construct a discrete collision mechanics on the points of the
grid in order to preserve the main physical properties.
Unfortunately DVM have the same computational cost of a product
quadrature rule and due to the particular choice of the nodes
imposed by the conservation properties the accuracy of the schemes
seems to be less than first order~\cite{BoPaSc:cons:97,
PaSc:stabcvgDVM:98, HePa:DVM:02}.

More recently a new class of methods based on the use of spectral
techniques in the velocity space has attracted the attention of
the scientific community. The method was first developed for kinetic
equations in~\cite{PePa:96}, inspired from spectral methods in
fluid mechanics~\cite{CHQ:88} and the use of Fourier transform
tools in the analysis of the Boltzmann
equation \cite{Boby:maxw:88}. It is based on a Fourier-Galerkin
approximation of the equation. Generalizations of the
method and spectral accuracy have been given
in~\cite{PaRu:spec:00, PaRu:stab:00}. This method, thanks to its
generality, has been applied also to non homogeneous
situations~\cite{FiRu:FBE:03}, to the Landau
equation~\cite{FiPa:02, PaRuTo:00} and to the case of granular
gases~\cite{NaGiPaTo:03, FiPa:03}. A related numerical strategy
based on the direct use of the fast Fourier transform
(FFT) has been developed in~\cite{BoRj:HS:97, BoRj:HS:99}.

The lack of discrete conservations in the spectral scheme (mass is
preserved, whereas momentum and energy are approximated with
spectral accuracy) is compensated by its higher accuracy and
efficiency. In fact it has been shown that these spectral schemes
permit to obtain spectrally accurate solutions with a reduction of
the computational cost strictly related to the particular
structure of the collision operator. A reduction from $O(N^2)$ to
$O(N\log_2 N)$ is readily deducible for the Landau equation,
whereas in the Boltzmann case such a reduction had been
obtained until now only at the price of a poor accuracy (in particular
the loss of the spectral accuracy), see~\cite{BoRj:HS:97, BoRj:HS:99}.

Finally we mention that spectral methods have been successfully
applied also to the study of {\em non cut-off} Boltzmann
equations, like for RGD in the grazing collision
limit~\cite{PaToVi:03} and for granular flows in the quasi-elastic
limit~\cite{NaGiPaTo:03}. In particular, during these asymptotic
processes it is possible to obtain intermediate approximations
that can be evaluated with fast algorithms that brings the overall
computational cost to $O(N \log_2 N)$. These idea has been used
in~\cite{Pa:03} to obtain fast approximated algorithms for the
Boltzmann equation.




For a recent introduction to numerical methods for the Boltzmann
equation and related kinetic equations we refer the reader
to~\cite{DPR:03}.

In this paper we shall focus on the two main questions in the
approximation the Boltzmann equation by deterministic schemes,
that is the computational complexity and the accuracy of the
numerical schemes for computing the collision operator $Q$.

Let us mention that a major problem associated with deterministic
methods that use a fixed discretization in the velocity domain is
that the velocity space is approximated by a finite region.
Physically the domain for the velocity is $\R^d$. But, as soon as
$d \ge 2$, the property of having compact support is not conserved
by the collision operator (in fact for some Boltzmann models in
dimension $d=1$, like granular models, the support is
conserved~\cite{NaGiPaTo:03}). In general the collision process ``spreads'' the support by a
factor $\sqrt{2}$ (see~\cite{PulvWenn:binf:97,M1:04}). As
a consequence, for the continuous equation in time, the function
$f$ is immediately positive in the whole velocity domain $\R^d$.

Thus at the numerical level some non physical condition has to be
imposed to keep the support of the function in velocity uniformly
bounded. In order to do this there are two main strategies, which we shall make more precise in the sequel.
\begin{enumerate}
\item One can remove the physical binary collisions that will lead outside the bounded
 velocity domain, which means a possible increase of the number of
 local invariants. If this is done properly ({\em i.e.} ``without removing too many collisions''),
 the scheme remains conservative (and without spurious invariants). However this truncation breaks
 down the convolution-like structure of the collision operator, which requires the invariance in velocity.
 Indeed the modified collision kernel depends on $v$ through the boundary conditions.
 This truncation is the starting point of most schemes based on discrete velocity models in a bounded domain.
\item One can add some non physical binary collisions by periodizing the function and
 the collision operator. This implies the loss of some local
 invariants (some non physical collisions are added).
 Thus the scheme is not conservative anymore, except for the mass
 if the periodization is done carefully (and possibly the momentum if some symmetry
 properties are satisfied by the function). In this way the structural properties of the collision operator
 are maintained and thus they can be exploited to derive fast algorithms. This periodization is the basis of
 the spectral method.
 \end{enumerate}

Note that in both cases by enlarging enough the computational domain
the number of removed or added collisions can be made
negligible (as it is usually done for removing the aliasing error of
the FFT, for instance see~\cite{CHQ:88}) as well as the error in the local invariants.

In this paper we shall focus on the second approach, which means
that the schemes have to deal with some aliasing error introduced
by the periodization. In this way, for a particular class of
interactions, using a Carleman-like representation of the collision
operator we are able to derive spectral methods
that can be evaluated through fast algorithms.
The class of interactions includes {\em Maxwellian molecules} in
dimension two and {\em hard spheres} molecules in dimension three.


The rest of the paper is organized in the following way. In
Section~\ref{sec:mp:rep} we introduce a Carleman-like representation of the
collision operator which is used as a starting point for the
development of our methods. After the derivation of the schemes
the details of the fast spectral algorithm together with its
accuracy properties are given in Section~\ref{sec:mp:fast}.
In a separate Appendix we show a possible way to extend
the present fast schemes to general collision interactions.

\section{Carleman-like representation and approximation of the collision operator}\label{sec:mp:rep}
\setcounter{equation}{0}

In this section we shall approximate the collision operator
starting from a representation which somehow conserves more
symmetries of the collision operator when one truncates it in a
bounded domain. This representation was used
in \cite{BoRj:HS:97,BoRj:HS:99,BoRj:maxw:98,IbRj:02} and it is
close to the classical Carleman representation
(cf.~\cite{Carl:fond:32}). Also the kind of periodization
inspired from this representation was implicitly used in \cite{BoRj:HS:99}.

\subsection{The Boltzmann collision operator in bounded domains}

The basic identity we shall need is
 \begin{equation}\label{MP:eq:dirac}
 \frac{1}{2} \, \int_{\ens{S}^{d-1}} F(|u|\sigma - u) \, d\sigma
 = \frac{1}{|u|^{d-2}} \, \int_{\R^d} \delta(2 \, x \cdot u + |x|^2) \, F(x) \, dx,
 \end{equation}
and can be verified easily by completing the square in the delta Dirac function, taking the spherical
coordinate $x=r \, \sigma$ and performing the change of variable $r^2 = s$.

\begin{figure}[h]
 \epsfysize=7cm
 $$\epsfbox{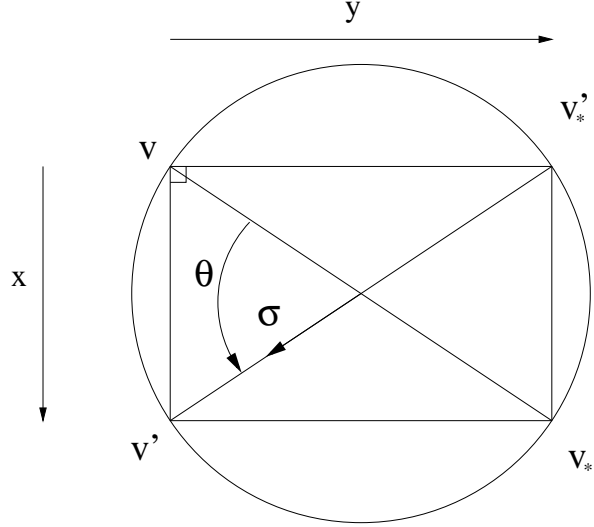}$$
 \caption{Geometry of the collision $(v,v_*)\leftrightarrow (v',v'_*)$.}\label{MP:fig:col}
 \end{figure}

Setting $u=v-v_*$ we can write the collision operator in the form
  \begin{multline*}
  Q (f,f)(v) = \int_{v_* \in \R^d}
  \Bigg\{ \int_{\sigma \in \ens{S}^{d-1}}  B(|u|,\cos \theta) \\
  \left[ f\Big(v_* - \frac{|u|\sigma - u}{2}\Big) \, f\Big(v+ \frac{|u|\sigma - u}{2}\Big)
  - f(v_*) \, f(v) \right] \, d\sigma \Bigg\} \, dv_*
  \end{multline*}
and thus equation~\eqref{MP:eq:dirac} yields
  \begin{multline*}
  Q (f,f)(v) = 2 \, \int_{v_* \in \R^d}
  \Bigg\{ \int_{x \in \R^d}  B\left(|u|,\frac{x \cdot u}{|x| |u|}\right) \, \frac{1}{|u|^{d-2}} \,
  \delta(2 \, x \cdot u + |x|^2) \\
  \Big[ f(v_* - x/2) \, f(v+ x/2) - f(v_*) \, f(v) \Big] \, dx \Bigg\} \, dv_*.
  \end{multline*}
Now let us make the change of variable $x \to x/2$ in $x$ to get
  \begin{multline*}
  Q (f,f)(v) = 2^{d+1} \, \int_{v_* \in \R^d}
  \int_{x \in \R^d}  B\left(|u|,\frac{x \cdot u}{|x| |u|}\right) \, \frac{1}{|u|^{d-2}} \,
  \delta(4 \, x \cdot u + 4 |x|^2) \\
  \left[ f(v_* - x) \, f(v+ x) - f(v_*) \, f(v) \right] \, dx \, dv_*
  \end{multline*}
and then setting $y=v_*-v-x$ in $v_*$ we obtain
  \begin{multline*}
  Q (f,f)(v) = 2^{d+1} \, \int_{y \in \R^d}
  \int_{x \in \R^d}  B\left(|u|,\frac{x \cdot u}{|x| |u|}\right) \, \frac{1}{|u|^{d-2}} \,
  \delta(- 4 x \cdot y) \\
  \left[ f(v + y) \, f(v+ x) - f(v+x+y) \, f(v) \right] \, dx \, dy
  \end{multline*}
where now $u=-(x+y)$. Thus in the end we have
  \begin{multline*}
  Q (f,f)(v) = 2^{d-1} \, \int_{x \in \R^d} \int_{y \in \R^d}
  B\left(|x+y|,-\frac{x \cdot (x+y)}{|x| |x+y|}\right) \, \frac{1}{|x+y|^{d-2}} \\
  \delta(x \cdot y) \,
  \left[ f(v + y) \, f(v+ x) - f(v+x+y) \, f(v) \right] \, dx \, dy.
  \end{multline*}
Figure~\ref{MP:fig:col} sums up the different geometrical quantities
of the usual representation and the one we derived from Carleman's
one.

Now let us consider the bounded domain $\mathcal{D}_T=[-T,T]^d$
($0< T < +\infty$). There are two possibilities of truncation to reduce the
collision process in a box. From now on let us write
 \[ \tilde{B}(x,y) = 2^{d-1} \, B\left(|x+y|,-\frac{x \cdot (x+y)}{|x| |x+y|}\right) \, |x+y|^{-(d-2)}. \]
One can easily see that on the manifold defined by $x \cdot y =0$,
a simpler formula is (using the parities of the collision kernel)
  \begin{multline}\label{eq:Btilde}
  \tilde{B}(x,y) = \tilde{B}(|x|,|y|) =
  2^{d-1} \, B\left(\sqrt{|x|^2+|y|^2}, \frac{|x|}{\sqrt{|x|^2+|y|^2}} \right) \, (|x|^2+|y|^2)^{-\frac{d-2}2}.
  \end{multline}
First one can remove the collisions connecting with some points
out of the box. This is the natural preliminary stage for deriving
conservative schemes based on the discretization of the velocity.
In this case there is no need for a truncation on the modulus of
$x$ and $y$ since we impose them to stay in the box.
It yields
  \begin{multline*}
  Q^{\mbox{tr}} (f,f)(v) = \int \int_{\big\{
  x, \, y \, \in \, \R^d \ | \ v+x, \, v+y, \, v+x+y \, \in \, \mathcal{D}_T \big\} }
  \tilde{B}(x,y) \,
  \delta(x \cdot y) \\
  \left[ f(v + y) \, f(v+ x) - f(v+x+y) \, f(v) \right] \, dx \, dy
  \end{multline*}
defined for $v \in  \mathcal{D}_T$. One can easily check that the
following weak form is satisfied by this operator
 \begin{multline}\label{eq:weaktr}
 \int Q^{\mbox{tr}} (f,f) \, \varphi(v) \, dv = \frac{1}{4} \,
 \int \int \int_{\big\{ v,\, x, \, y \, \in \, \R^d \ | \
 v,\, v+x,\, v+y,\, v+x+y \, \in \, \mathcal{D}_T \big\}}
 \tilde{B}(x,y) \, \delta(x \cdot y) \\
 f(v+x+y) \, f(v) \left[ \varphi(v + y) + \varphi(v+ x) - \varphi(v+x+y) - \varphi(v) \right] \, dv \, dx \, dy
 \end{multline}
and this implies conservation of mass, momentum and energy as well
as the $H$ theorem on the entropy. Note that at this level this
formulation gives no advantage with respect to the usual one
obtained from~\eqref{eq:mp:Q} by restricting $v, v_*, v', v'_* \in
\mathcal{D}_T$ (except that consistency results for discrete
velocity models seem easier to prove when they are derived by
quadrature on this formulation, see~\cite{HePa:DVM:02}). The
problem of this truncation on a bounded domain is the fact that we
have changed the collision kernel itself by adding some artificial
dependence on $v, v_*,v',v'_*$. In this way convolution-like
properties are broken.

A different approach consists in periodizing the function $f$ on
the domain $\mathcal{D}_T$. This amounts in adding some
non-physical collisions by connecting some points in the domain
$\mathcal{D}_T$ which are geometrically included in a collision
circle ``modulo $T$'' ({\em i.e.} up to a translation of $T$ of
certain points in certain directions). Here we have to truncate
the integration in $x$ and $y$ since periodization would yield
infinite result if not. Thus we set them to vary in
$\mathcal{B}_R$, the ball of center $0$ and radius $R$. For a
compactly supported function $f$ with support $\mathcal{B}_S$, we
take $R={2S}$ in order to obtain all possible collisions.
Then a geometrical argument (see~\cite{PaRu:spec:00}) shows that
using the periodicity of the function it is enough to take $T \ge
(1 +3\sqrt{2})S/2$ to prevent intersections of the regions where
$f$ is different from zero. Note that here this so-called
dealiasing condition is slightly worst from the one in
\cite{PaRu:spec:00}, since the truncation on the modulus of $x$
and $y$ in the ball $\mathcal{B}_R$ implies only a truncation in
the ball $\mathcal{B}_{\sqrt{2}R}$ for the relative velocity.

The operator now reads
 \begin{multline}\label{eq:repR}
 Q^R (f,f)(v) = \int_{x \in \mathcal{B}_R} \int_{y \in \mathcal{B}_R}
 \tilde{B}(x,y) \, \delta(x \cdot y) \\
 \left[ f(v + y) f(v+ x) - f(v+x+y) f(v) \right] \, dx \, dy
 \end{multline}
for $v \in \mathcal{D}_T$ (the expression for $v \in \R^d$ is
deduced by periodization). The interest of this representation is
to preserve the real collision kernel and its properties.

By making some translation changes of variable on $v$ (by $x$, $y$
and $x+y$), using the changes $x \to -x$ and $y \to -y$ and the
fact that
 \[ \tilde{B}(-x,y) \, \delta(-x \cdot y) = \tilde{B}(x,y) \, \delta(x \cdot y)
 = \tilde{B}(x,-y) \, \delta(x \cdot -y) \]
one can easily prove that for any function $\varphi$ {\em
periodic} on $\mathcal{D}_T$ the following weak form is satisfied
 \begin{multline}\label{eq:weakpe}
 \int_{\mathcal{D}_T} Q^R (f,f) \, \varphi(v) \, dv = \frac{1}{4} \,
 \int_{v \in \mathcal{D}_T} \int_{x \in \mathcal{B}_R} \int_{y \in \mathcal{B}_R}
 \tilde{B}(x,y) \, \delta(x \cdot y) \\
 f(v+x+y) f(v) \left[ \varphi(v + y) + \varphi(v+ x) - \varphi(v+x+y) - \varphi(v) \right] \, dv \, dx \,
 dy.
 \end{multline}

About the conservation properties one can shows that
 \begin{enumerate}
 \item The only invariant $\varphi$ is $1$: it is the only periodic function on $\mathcal{D}_T$
 such that
  \[ \varphi(v + y) + \varphi(v+ x) - \varphi(v+x+y) - \varphi(v) = 0 \]
 for any $v \in \mathcal{D}_T$ and $x \bot y \in \mathcal{B}_R$ (see~\cite{Cerc:75} for instance).
 It means that the mass is locally conserved but not necessarily the momentum and energy.
 \item When $f$ is even there is {\em global} conservation of momentum, which is $0$ in this case.
 Indeed $Q^R$ preserves the parity property of the solution, which can be checked using
 the change of variable $x \to -x$, $y \to -y$.
 \item The collision operator satisfies formally the $H$ theorem
   \[ \int_{\R^d} Q^R (f,f)\log(f) \, dv \leq 0. \]
 \item If $f$ has compact support included in $\mathcal{B}_S$, and we have
 $R =  2S$ and $T \ge (3\sqrt{2}+1)S/2$ (no aliasing condition, see \cite{PaRu:spec:00}
 for a detailed discussion), then no unphysical collisions occur
 and thus mass, momentum and energy are preserved. Obviously this compactness is not preserved with
 time since the collision operator spreads the support of $f$ by a factor $\sqrt{2}$.
 \end{enumerate}
To sum up one could say that the lack of conservations originates
from the fact that the geometry of the collision does not respect
the periodization.

Finally we give the Cauchy theorems for the homogeneous Boltzmann
equations in $\mathcal{D}_T$ computed with $Q^{\mbox{tr}}$ or
$Q^R$.
 \begin{theorem}\label{theo:cauchy}
 Let $f_0 \in L^1 (\mathcal{D}_T)$ be a nonnegative function. Then there exists a unique
 solution $f \in C^1(\R_+, L^1 (\mathcal{D}_T))$ to the Cauchy problems
  \begin{equation}\label{eq:cauchyTR}
  \derpar{f}{t} = Q^{\mbox{{\em tr}}} (f,f), \ \ \ f(t=0, \cdot)= f_0
  \end{equation}
  \begin{equation} \label{eq:cauchyPER}
  \derpar{f}{t} = Q^R (f,f), \ \ \ f(t=0, \cdot)= f_0
  \end{equation}
 which is nonnegative and has constant mass (and so constant $L^1$ norm). If $f_0$ has finite entropy,
 the entropy is finite and non-decreasing for all time.
 Moreover in the case~\eqref{eq:cauchyTR},
 if $f_0$ has finite momentum (respectively energy) on $\mathcal{D}_T$, the momentum
 (respectively energy) is conserved with time.
 \end{theorem}

\Remark When the initial data $f_0$ is nonnegative and has finite
mass and entropy, it is possible to show by the Dunford-Pettis
compactness theorem that the solution $f$ converges weakly in
$L^1(\mathcal{D}_T)$, as $t$ goes to infinity, to the unique
maximum of the entropy functional compatible with the conservation
law(s) (and the periodicity in the case~\eqref{eq:cauchyPER}). In
the case~\eqref{eq:cauchyTR} this equilibrium state is a sort of
truncated Maxwellian on $\mathcal{D}_T$ defined by the
conservation laws (see~\cite{Cerc:75}). In the
case~\eqref{eq:cauchyPER} this equilibrium state is a constant
defined by the mass of the initial data, which is due to the
effect of aliasing in the very long-time. We omit the proof for
brevity.

\begin{proof}[Proof of Theorem~\ref{theo:cauchy}]
For clarity we briefly sketch the main lines of the proof. The
existence and uniqueness are proved by the method of Arkeryd for
bounded collision kernels, see~\cite[Part~I,
Proposition~1.1]{Arke:I+II:72}. In our case the collision kernel
is bounded because of the boundedness of the domain. The only a
priori estimate required in~\cite[Part~I,
Proposition~1.1]{Arke:I+II:72} is the mass conservation, valid for
the two equations under consideration. This method is based on a
monotonicity argument to prove propagation of the sign of the
solution. The argument relies on a splitting of the collision
operator $Q$ into a gain part $Q^+$ which is monotonic ({\em i.e.} $Q^+
(f,f)$ is non-negative when $f$ is non-negative), and a loss part
$Q^-$ which writes $Q^- (f,f) = L(f) f$ with $L$ is a linear
operator such that $\| L(f) \|_{\infty} \le C \, \| f \|_{L^1}$.
One can check easily that this splitting is still valid for the
two collision operators $Q^{\mbox{tr}}$ and $Q^R$. For brevity we
omit the details and refer to the article~\cite{Arke:I+II:72}. The
conservation law(s) and the $H$ theorem are deduced from the weak
forms~\eqref{eq:weaktr} and~\eqref{eq:weakpe} (see the proof
of~\cite[Part~I, Proposition~1.2]{Arke:I+II:72} and~\cite[Part~I,
Theorem~2.1]{Arke:I+II:72}).
\end{proof}




\subsection{Application to spectral methods}

In this Section we use the representation $Q^R$ to derive new
spectral methods. The spectral methods for kinetic equations
originated in the works of~\cite{PePa:96,PaRu:spec:00}, and were
further developed in~\cite{PaRu:stab:00,FiRu:FBE:03}. Before they
had a long history in fluid mechanics, see~\cite{CHQ:88}.

The main change compared to the usual spectral method is in the
way we truncate the collision operator. In fact as we shall see in
the next section this yields better decoupling properties between
the arguments of the operator.

To simplify notations let us take $T=\pi$. Hereafter we use just
one index to denote the $d$-dimensional sums of integers.

The approximate function $f_N$ is represented as the truncated Fourier series
  \[ \left\{
     \begin{array}{l}\displaystyle
     f_N (v) = \sum_{k=-N}^{N} \hat{f}_k e^{i k \cdot v}, \vspace{0.2cm} \\ \displaystyle
     \hat{f}_k = \frac{1}{(2 \pi)^d} \, \int _{\mathcal{D}_\pi} f(v) e^{-ik \cdot v} \, dv.
     \end{array}
     \right. \]
The spectral equation is the projection of the collision equation in $\ens{P}^N$, the $(2N+1)^d$-dimensional
vector space of trigonometric polynomials of degree at most $N$ in each direction, {\em i.e.}
 \begin{equation*}
 \derpar{f_N}{t} = \mathcal{P}_N Q^R (f_N,f_N)
 \end{equation*}
where $\mathcal{P}_N$ denotes the orthogonal projection on
$\ens{P}^N$ in $L^2(\mathcal{D}_\pi)$.
A straightforward computation leads
to the following set of ordinary differential equations on the
Fourier coefficients
 \begin{equation}\label{eq:ode}
 \hat{f}_k '(t) =
 \sum_{\underset{l+m=k}{l,m=-N}}^{N} \hat{\beta}(l,m) \, \hat{f}_l \, \hat{f}_m, \ \ \
 k=-N,...,N
 \end{equation}
where $\hat{\beta}(l,m)$ are the so-called {\em kernel modes}, given by
 \begin{equation*}
 \hat{\beta}(l,m) = \int_{x \in \mathcal{B}_R} \int_{y \in \mathcal{B}_R}
 \tilde{B}(x,y) \, \delta(x \cdot y) \, \left[
 e^{i l \cdot x} \, e^{i m \cdot y} - e^{i m \cdot (x+y)} \right] \, dx \, dy.
 \end{equation*}
The kernel modes can be written as
 \begin{equation*}
 \hat{\beta}(l,m) = \beta(l,m) - \beta(m,m)
 \end{equation*}
where
 \begin{equation*}
 \beta(l,m) = \int_{x \in \mathcal{B}_R} \int_{y \in \mathcal{B}_R}
 \tilde{B}(x,y) \, \delta(x \cdot y) \,
 e^{i l \cdot x} \, e^{i m \cdot y} \, dx \, dy.
 \end{equation*}
Therefore in the sequel we shall focus on $\beta$, and one easily
checks that $\beta(l,m)$ depends only on $|l|$, $|m|$ and $|l \cdot
m|$.

Note that the usual way to truncate the Boltzmann collision
operator for periodic function starts from the following
representation (see \cite{PaRu:spec:00})
 \begin{multline}\label{eq:repusual}
 Q (f,f) = \int_{u \in \R^d} \int_{\sigma \in \ens{S}^{d-1}}
           B(|u|,\cos \theta) \\
           \Big[ f\big(v- (u - |u|\sigma)/2 \big) f\big(v- (u + |u|\sigma)/2 \big)
           - f(v) f( v -u) \Big] \, d\sigma \, du
 \end{multline}
and then truncate the parameter $u=x+y$ in order that $u \in
\mathcal{B}_{R}$. Thus we have
 \begin{multline*}
 Q^R _{\mbox{\scriptsize{usual}}} (f,f)(v) = \int_{x \in \R^d} \int_{y \in \R^d}
 \tilde{B}(x,y) \, \delta(x \cdot y) \, {\bf \chi}_{\{|x+y|\le R\}}  \\
 \left[ f(v + y) f(v+ x) - f(v+x+y) f(v) \right] \, dx \, dy
 \end{multline*}
where ${\bf \chi}_{\{|x+y|\le R\}}$ denotes the characteristic
function of the set $\{|x+y|\le R\}$.
One can notice that here $x$ and $y$ are also restricted to the
ball $\mathcal{B}_R$ but the condition $|x+y|^2 = |x|^2 + |y|^2
\le R^2$ couples the two modulus, such that the ball is not
completely covered (for instance, if $x$ and $y$ have both modulus
$R$, the condition is not satisfied, since $|x+y| = \sqrt{2} R$).

Finally let us compare the new kernel modes with the usual ones.
As a consequence of the representation~\eqref{eq:repusual}, the
usual kernel modes (cf.~\cite{PaRu:spec:00}) are
  \begin{multline*}
  \hat{\beta} _{\mbox{\scriptsize{usual}}} (l,m) = \int_{u \in \mathcal{B}_R} \int_{\sigma \in \ens{S}^{d-1}}
            B(|u|,\cos \theta) \,
            \Big[ e^{-i\frac{u \cdot (l+m) + |u| \sigma \cdot (m-l)}{2}}
            - e^{-i (u \cdot m)}  \Big] \, d\sigma \, du
  \end{multline*}
and hence coming back to the representation in $x$ and $y$,
  \begin{multline*}
  \hat{\beta} _{\mbox{\scriptsize{usual}}} (l,m) = \int_{x \in \mathcal{B}_R} \int_{y \in \mathcal{B}_R}
  \tilde{B}(x,y) \, \delta(x \cdot y) \, {\bf \chi}_{\{|x+y| \le R\}} \,
  \left[ e^{i l \cdot x} \, e^{i m \cdot y} - e^{i m \cdot (x+y)} \right] \, dx \, dy.
  \end{multline*}
Thus the usual representation contains more coupling between $x$
and $y$ and it is less appropriate for the construction of fast
algorithms.

\section{Fast spectral algorithm for a class of collision kernels}\label{sec:mp:fast}
\setcounter{equation}{0}

As soon as one is searching for fast deterministic algorithms for
the collision operator, {\em i.e.} algorithm with a cost lower than
$O(N^{2d+\var})$ (which is the cost of a usual discrete velocity
model, with typically $\var=1$), one has to find some way to
compute the collision operator {\em without going through all the
couples of collision points} during the computation. This leads
naturally to search for some convolution structure (discrete or
continuous) in the operator. Unfortunately, as discussed in the
previous sections, this is rather contradictory with the search
for a conservative scheme in a bounded domain, since the boundary
condition needed to prevent for the outgoing or ingoing collisions
breaks the invariance. Thus fast algorithms seem more adapted to
spectral methods, or more in general to methods where the
invariance is conserved thanks to the periodization.

Here we search for a convolution structure in the equations~\eqref{eq:ode}.
The aim is to approximate each $\hat{\beta}(l,m)$ by a sum
 \[ \hat{\beta}(l,m) \simeq \sum_{p=1} ^{A} \alpha_p (l) \alpha' _p (m). \]
This gives a sum of $A$ discrete convolutions and so the algorithm
can be computed in $O(A \, N^d \log_2 N)$ operations by means of
standard FFT techniques~\cite{CHQ:88, CoTu:65}. Obviously this is
equivalent to obtain such a decomposition on $\beta$. To this
purpose we shall use a further approximated collision operator
where the number of possible directions of collision is reduced to
a finite set.

The starting point of our study is an idea of \cite{BoRj:HS:99}: use the
Carleman-like representation \eqref{eq:repR} to obtain a convolution structure
for every fixed directions of the vectors $x$ and $y$. In this work \cite{BoRj:HS:99}
the corresponding set of directions
  \[ S= \left\{ (e,e') \in \ens{S}^{N-1} \times \ens{S}^{N-1} \ | \
              e \bot e' \right\} \]
is very difficult to discretize in a way that preserves the symmetry properties
of the collision operator.
No systematic process is available and the discretization
is done only for some particular number of grid points. Then the FFT is used in each couple of direction and
finally a correction is imposed at the end to preserve the conservation laws. However no
consistency result is available and the accuracy suggested by the numerical simulations
is of order $1$. The two main new ingredients of our method are:
\begin{itemize}
\item First we project the collision operator on the Fourier basis. This enables to integrate one
of the two coordinates of the manifold $S$ and to reduce to the discretization of the
sphere $\ens{S}^{N-1}$. This discretization is straightforward and can be made easily
to preserve the symmetries of the collision operator. Moreover it reduces the complexity
of the algorithm by suppressing $N-2$ degrees of freedom to discretize.
\item Second we choose to discretize $\ens{S}^{N-1}$ by the rectangular rule. Indeed the
periodization shall imply that this quadrature rule is of infinite order. This point
will allow to obtain a spectrally accurate scheme, and adaptativity properties.
\end{itemize}

\subsection{A semi-discrete collision operator}

We write $x$ and $y$ in spherical coordinates
 \begin{multline}\label{eq:sansapprox}
 Q^R (f,f)(v) = \frac14 \, \int_{e \in \ens{S}^{d-1}} \int_{e' \in \ens{S}^{d-1}}
 \, \delta(e \cdot e') \, de \, de' \\
 \Bigg\{ \int_{-R} ^R \int_{-R} ^R \rho^{d-2} \, (\rho')^{d-2} \,
 \tilde{B}(\rho, \rho')
 \big[ f(v + \rho' e') f(v+ \rho e) - f(v+\rho e + \rho' e') f(v) \big]
 \, d\rho \, d\rho' \Bigg\}.
 \end{multline}
Let us take $\mathcal{A}$ a set of orthogonal couples of unit
vectors $(e,e')$, which is even: $(e,e') \in \mathcal{A}$ implies
that $(-e,e')$, $(e,-e')$ and $(-e,-e')$ belong to $\mathcal{A}$ (this property on the
set $\mathcal{A}$ is required to preserve the conservation properties of the operator).
Now we define $Q_R ^\mathcal{A}$ to be
 \begin{multline*}
 Q^{R,\mathcal{A}} (f,f)(v) = \frac14 \, \int_{(e,e') \in \mathcal{A}}
 \Bigg\{ \int_{-R} ^R \int_{-R} ^R \rho^{d-2} \, (\rho')^{d-2} \,
 \tilde{B}(\rho, \rho') \\ \big[ f(v + \rho' e') f(v+ \rho e) -
 f(v+\rho e + \rho' e') f(v) \big]
 \, d\rho \, d\rho' \Bigg\} \, d\mathcal{A}
 \end{multline*}
where $d\mathcal{A}$ denotes a measure on $\mathcal{A}$
which is also even in the sense that $d\mathcal{A}(e,e')=
d\mathcal{A}(-e,e')=d\mathcal{A}(e,-e')= d\mathcal{A}(-e,-e')$.
Using translation changes of variable on $v$ by $\rho e$,
$\rho' e'$ and $\rho e + \rho' e'$ and the symmetries of the set
$\mathcal{A}$ one can easily derive the following weak form on
$Q_R^\mathcal{A}$. For any function $\varphi$ {\em periodic} on
$\mathcal{D}_T$,
 \begin{multline*}
 \int_{\mathcal{D}_T} Q^{R,\mathcal{A}} (f,f) \, \varphi(v) \, dv  =
 \frac1{16} \, \int_{v \in \mathcal{D}_T} \int_{(e,e') \in \mathcal{A}}
 \int_{-R} ^R \int_{-R} ^R \rho^{d-2} \, (\rho')^{d-2} \,
 \tilde{B}(\rho, \rho') \\
 f(v+\rho e +\rho' e') f(v) \, \Big[ \varphi(v +\rho' e') + \varphi(v+ \rho e ) -
 \varphi(v+\rho e+\rho' e') - \varphi(v) \Big] \, d\rho \, d\rho' \,
 d\mathcal{A} \, dv.
 \end{multline*}
This immediately gives the same conservations properties as $Q^R$.
Of course one could also prove exactly as for $Q^R$:
 \begin{theorem}
 Let $f_0 \in L^1 (\mathcal{D}_T)$ be a nonnegative function. Then there exists a unique
 solution $f \in C^1(\R_+, L^1 (\mathcal{D}_T))$ to the Cauchy problem
  \[ \derpar{f}{t} = Q^{R,\mathcal{A}} (f,f), \ \ \ f(t=0, \cdot)= f_0  \]
 which is nonnegative and has constant mass (and so constant $L^1$ norm).
 Moreover, if $f_0$ has finite entropy, the entropy is non-decreasing with time.
 \end{theorem}



\subsection{Expansion of the kernel modes}

We make the {\em decoupling assumption} that
 \begin{equation}\label{eq:decoup}
 \tilde{B}(x,y) = a(|x|) \, b(|y|).
 \end{equation}


This assumption is obviously satisfied if $\tilde B$ is constant.
This is the case of Maxwellian molecules in dimension two, and hard
spheres in dimension three (the most relevant kernel for
applications). Extensions to more general interactions are
discussed in the Appendix.

First let us deal with dimension $2$
with $\tilde{B} =1$ to explain the method. Here we write $x$ and
$y$ in spherical coordinates $x = \rho e$ and $y = \rho' e'$ to
get
 \begin{equation*}
 \beta (l,m) = \frac14 \, \int_{e \in \ens{S}^1} \int_{e' \in \ens{S}^1}
 \delta(e \cdot e') \,
 \left[ \int_{-R} ^R e^{i \rho (l \cdot e)} \, d\rho \right] \,
 \left[ \int_{-R} ^R e^{i \rho' (m \cdot e')} \, d\rho' \right] \, de \,
 de'.
 \end{equation*}
Let us denote by
 \[ \phi_R ^2 (s) = \int_{-R} ^R e^{i \rho s} \, d\rho, \]
for $s \in \R$. It is easy to see that $\phi_R ^2$ is even and we can give the explicit
formula
 \[ \phi_R ^2 (s) = 2 \, R \, \mbox{Sinc} (R s) \]
with $\mbox{Sinc}(\theta) = (\sin \theta)/\theta$.

Thus we have
 \begin{equation*}
 \beta (l,m) = \frac14 \, \int_{e \in \ens{S}^1} \int_{e' \in \ens{S}^1}
 \delta(e \cdot e') \, \phi_R ^2 (l \cdot e) \, \phi_R ^2 (m \cdot e') \, de \, de'
 \end{equation*}
and thanks to the parity property of $\phi_R ^2$ we can adopt the following periodic
parametrization
 \begin{equation*}
 \beta (l,m) =  \int_0 ^{\pi} \phi_R ^2 (l \cdot e_{\theta})\, \phi_R ^2 (m \cdot e_{\theta+\pi/2}) \,
 d\theta.
 \end{equation*}
The function $\theta \to \phi_R ^2 (l \cdot e_{\theta})\, \phi_R
^2 (m \cdot e_{\theta+\pi/2})$ is periodic on $[0,\pi]$ and thus
the rectangular quadrature rule is of infinite order and optimal.
A regular discretization of $M$ equally spaced points thus gives
 \begin{equation*}
 \beta (l,m) = \frac{\pi}{M} \, \sum_{p=0} ^{M-1} \alpha_p (l) \alpha' _p (m)
 \end{equation*}
with
 \[ \alpha _p (l) = \phi_R ^2 (l \cdot e_{\theta_p}), \hspace{0.8cm}
    \alpha' _p (m) = \phi_R ^2 (m \cdot e_{\theta_p+\pi/2}) \]
where $\theta_p = \pi p/M$.

More generally under the decoupling assumption~\eqref{eq:decoup}
on $\tilde{B}$, we get the following decomposition formula
 \begin{equation*}
 \beta (l,m) = \frac{\pi}{M} \, \sum_{p=0} ^{M-1} \alpha_p (l) \alpha' _p (m)
 \end{equation*}
where
 \[ \alpha_p (l) = \phi_{R,a} ^2 (l \cdot e_{\theta_p}), \hspace{0.8cm}
    \alpha ' _p (m) = \phi_{R,b} ^2 (m \cdot e_{\theta_p+\pi/2}) \]
and
 \[ \phi_{R,a} ^2 (s) = \int_{-R} ^R a(\rho) \,  e^{i \rho s} \, d\rho , \hspace{0.8cm}
    \phi_{R,b} ^2 (s) = \int_{-R} ^R b(\rho') \,  e^{i \rho' s} \, d\rho' \]
with $\theta_p = \pi p/M$.

\Remark In the symmetric case $a=b$ (for instance for hard spheres)
it is possible to parametrize $\beta(l,m)$ as
 \[  \beta (l,m) =  2 \, \int_0 ^{\pi/2} \phi_{R,a} ^2 (l \cdot e_{\theta})\,
     \phi_{R,a} ^2 (m \cdot e_{\theta+\pi/2}) \, d\theta \]
and the function $\theta \to \phi_{R,a} ^2 (l \cdot e_{\theta})\, \phi_{R,a} ^2 (m \cdot e_{\theta+\pi/2})$
is periodic on $[0,\pi/2]$.
Thus the decomposition can be obtained by applying the rectangular rule on this interval. At the
numerical level it yields a reduction of the cost by a factor $2$.

Now let us deal with dimension $d=3$ with $\tilde{B}$ satisfying the decoupling assumption~\eqref{eq:decoup}.
First we change to the spherical coordinates
 \begin{multline*}
 \beta (l,m) = \frac14 \, \int_{e \in \ens{S}^2} \int_{e' \in \ens{S}^2}
 \delta(e \cdot e') \,
 \left[ \int_{-R} ^R |\rho|  \, a(\rho) \, e^{i \rho (l \cdot e)} \, d\rho \right] \,
 \left[ \int_{-R} ^R |\rho'| \, b(\rho') \, e^{i \rho' (m \cdot e')} \, d\rho' \right] \, de \, de'
 \end{multline*}
and then we integrate first $e'$ on the intersection of the unit
sphere with the plane $e^\bot$,
 \begin{equation*}
 \beta (l,m) = \frac14 \, \int_{e \in \ens{S}^2} \phi_{R,a} ^3 (l \cdot e) \, \left[
 \int_{e' \in \ens{S}^2 \cap e^\bot} \phi_{R,b} ^3 (m \cdot e') \, de' \right] \, de
 \end{equation*}
where
 \begin{equation*}
 \phi_{R,a} ^3 (s) = \int_{-R} ^R |\rho|  \, a(\rho) \, e^{i \rho s} \, d\rho, \hspace{1cm}
 \phi_{R,b} ^3 (s) = \int_{-R} ^R |\rho|  \, b(\rho) \, e^{i \rho s} \, d\rho.
 \end{equation*}
Thus we get the following decoupling formula with two degrees of freedom
 \begin{equation*}
 \beta (l,m) = \int_{e \in \ens{S}^2 _+} \phi_{R,a} ^3 (l \cdot e) \,
 \psi_{R,b} ^3 \big(\Pi_{e^\bot}(m)\big) \, de
 \end{equation*}
where $\ens{S}^2 _+$ denotes the half-sphere and
 \begin{equation*}
 \psi_{R,b} ^3 \big(\Pi_{e^\bot}(m)\big) = \int_0 ^\pi \phi_{R,b} \big(|\Pi_{e^\bot}(m)|
 \, \cos\theta \big)  \, d\theta,
 \end{equation*}
(this formula can be derived performing the change of variable $de' = \sin \theta \, d\theta \, d\varphi$
with the basis $(e,u=\Pi_{e^\bot}(m)/|\Pi_{e^\bot}(m)|,e \times u)$).

Again in the particular case where $\tilde{B}=1$ (hard spheres
model), we can compute explicitly the functions $\phi_R ^3$ (in
this case $a=b=1$),
 \begin{equation*}
 \phi_R ^3 (s) = R^2 \, \left[ 2 \mbox{Sinc} (R s) - \mbox{Sinc} ^2 (R s /2) \right]. 
 \end{equation*}

Now the function $e \to \phi_{R,a} ^3 (l \cdot e) \, \psi_{R,b} ^3 \big(\Pi_{e^\bot}(m)\big)$ is periodic
on $\ens{S}^2 _+$ and so the rectangular rule is of infinite order and optimal.
Taking a spherical parametrization $(\theta,\varphi)$
of $e \in \ens{S}^2 _+$ and uniform grids of respective size $M_1$ and $M_2$ for $\theta$ and $\varphi$ we get
 \begin{equation*}
 \beta (l,m) = \frac{\pi^2}{M_1 M_2} \, \sum_{p,q=0} ^{M_1,M_2} \alpha_{p,q} (l) \alpha' _{p,q} (m)
 \end{equation*}
where
 \[ \alpha_{p,q} (l) = \phi_{R,a} ^3 \big(l \cdot e_{(\theta_p,\varphi_q)}\big),
 \hspace{0.8cm} \alpha' _{p,q} (m) = \psi_{R,b} ^3 \Big(\Pi_{e_{(\theta_p,\varphi_q)} ^\bot} (m)\Big) \]
and
 \[ (\theta_p,\varphi_q) = \Big(\frac{p \, \pi}{M_1}, \frac{q \, \pi}{M_2} \Big). \]
From now on we shall consider this expansion with $M=M_1=M_2$ to
avoid anisotropy in the computational grid.

\Remarks

1. It is possible to give more general exact formula in dimension $2$ and $3$
when $a(r) = |r|^t$, $b(r) = |r|^{t'}$ with $t,t' \in \N$
by computing derivatives along along $s$ of the two quantities
 \begin{equation*}
 \int_0 ^R \sin (\rho s) \, d\rho, \qquad
 \int_0 ^R \cos (\rho s) \, d\rho.
 \end{equation*}
\smallskip

2. For any dimension, we can construct as above an approximated collision operator $Q^{R,
\mathcal{A}_M}$ with
 \begin{equation*}
 \mathcal{A}_M = \Big\{ (e,e') \in \ens{S}^{d-1} \times \ens{S}^{d-1} \ \big| \
                        e \in \ens{S}^{d-1} _{M,+}, \ \ e' \in e^\bot \cap \ens{S}^{d-1} \Big\}
 \end{equation*}
where $\ens{S}^{d-1} _{M,+}$ denotes a uniform angular
discretization of the half sphere with $M$ points in each angular
coordinate (the other half sphere is obtained by parity). Let us
remark that this discretization contains exactly $M^{d-1}$ points.
From now on we shall denote
 \begin{equation*}
 Q^{R,M} = Q^{R, \mathcal{A}_M} = \sum_{p=1} ^{M^{d-1}} Q^{R,M} _p.
 \end{equation*}

\subsection{Spectral accuracy}

In this paragraph we are interested in computing the accuracy of
the scheme according to the three parameters $N$ (the number of
modes), $R$ (the truncation parameter), and $M$ (the number of
angular directions for each angular coordinate). Instead of
looking at the error on each kernel mode it is more convenient to
look at the error on the global operator. Here the Lebesgue spaces
$L^p$, $p=1 \ldots +\infty$, and the periodic Sobolev spaces $H^k
_p$, $k=0 \ldots +\infty$ refer to $\mathcal{D}_\pi$.

In order to give a consistency result, the first step will be
to prove a consistency result for the approximation of $Q^R$ by $Q^{R,M}$.
 \begin{lemma}\label{lem:rect}
 The error on the approximation of the collision operator is spectrally small, {\em i.e.} 
 for all $k > d-1$ such that $f \in H^k _p$
  \begin{equation*}
  \| Q^R (g,f) -  Q^{R,M} (g,f) \|_{L^2} \le C_1 \, \frac{ R^k \|g\| _{H^k _p} \|f\| _{H^k _p} }{M^k}.
  \end{equation*}
 \end{lemma}
\begin{proof}[Proof of Lemma~\ref{lem:rect}]
Starting from~\eqref{eq:sansapprox}, one gets
 \begin{multline*}
 Q^R (g,f)(v) = \frac12 \, \int_{e \in \ens{S}^{d-1} _+} \Bigg[ \int_{e' \in \ens{S}^{d-1} \cap e^\bot}
 \int_{-R} ^R \int_{-R} ^R \rho^{d-2} \, (\rho')^{d-2} \,
 \tilde{B}(\rho, \rho') \\
 \left[ g(v + \rho' e') f(v+ \rho e) - g(v+\rho e + \rho' e') f(v) \right]
 \, d\rho \, d\rho'  \, de' \Bigg] \, de.
 \end{multline*}
As the function in the brackets is a periodic function of $e$ on $\ens{S}^{d-1} _+$ with
period $\pi$ in each coordinate, one can apply the error estimate for the rectangular rule
(see for instance~\cite[Theorem~19.10]{Schat:91}). This error estimate is valid for $k > d-1$
and depends on the derivative along $e$ of this functional on the following way
 \begin{multline*}
 \| Q^R (g,f) -  Q^{R,M} (g,f) \|_{L^2} \le \frac{C}{2^k M^k} \, \sum_{i=1} ^{d-1}
 \Bigg\| \int_{e \in \ens{S}^{d-1} _+} \Bigg| \partial^k _{e_i} \int_{e' \in \ens{S}^{d-1} \cap e^\bot}
 \int_{-R} ^R \int_{-R} ^R \rho^{d-2} \, (\rho')^{d-2} \\
 \tilde{B}(\rho, \rho') \,
 \left[ g(v + \rho' e') f(v+ \rho e) - g(v+\rho e + \rho' e') f(v)
 \right]
 \, d\rho \, d\rho'  \, de' \Bigg| \, de \Bigg\|_{L^2 _v}
 \end{multline*}
where the constant is independent on $k$ and $\partial^k _{e_i}$
is the derivative of order $k$ along the coordinate $e_i$. Then a
straightforward computation gives
 \begin{multline*}
 \| Q^R (g,f) -  Q^{R,M} (g,f) \|_{L^2}  \le \frac{C R^k}{2^k M^k} \, \sum_{i=1} ^{d-1}
 \Bigg[ \sum_{k'+k''=k} \left( \begin{smallmatrix}  k \\ k'  \end{smallmatrix} \right)
 \big( \|Q^{R,+} (|\partial^{k'} g|, |\partial^{k''} f|) \|_{L^2} \\
 + \| Q^{R,-} (|\partial^{k'} g|, |\partial^{k''} f|)\|_{L^2}  \big) \Bigg]
 \end{multline*}
where $\partial^{k'}$ and $\partial^{k''}$ denote some derivatives of order $k'$ and $k''$.
Then using the estimates
 \[ \|Q^{R,+} (g,f), \ Q^{R,-} (g,f) \|_{L^2} \le C \, \|g\|_{L^2} \|f\|_{L^2} \]
proved in~\cite{FiRu:pre}\footnote{Which are consequences of the $L^p$ estimates
proved in~\cite{Gust:L^p:86,Gust:L^p:88}, and revisited in~\cite{MV:04}.}, we get
 \begin{equation*}
 \| Q^R (g,f) -  Q^{R,M} (g,f) \|_{L^2} \le \frac{C R^k}{M^k} \, \|g\| _{H^k _p} \|f\| _{H^k _p}
 \end{equation*}
which concludes the proof.
\end{proof}

For the second step we shall use the consistency result~\cite[Corollary~5.4]{PaRu:spec:00} on the operator $Q^{R}$,
which we quote here for the sake of clarity.
 \begin{lemma}\label{lem:PaRu}
 For all $k \in \N$ such that $f \in H^k _p$,
   \begin{equation*}
  \| Q^{R} (f,f) -  \mathcal{P}_N Q^{R} (f_N,f_N) \|_{L^2} \le \frac{C_2}{N^k}
                                     \, \Big( \|f\| _{H^k _p} + \|Q^R (f_N,f_N)\|_{H^k _p} \Big).
  \end{equation*}
 \end{lemma}
Combining these two results, one gets the following consistency result
 \begin{theorem}\label{theo:consist}
 For all $k > d-1$ such that $f \in H^k _p (\mathcal{D}_\pi)$,
  \begin{multline*}
  \| Q^R (f,f) -  \mathcal{P}_N Q^{R,M} (f_N,f_N) \|_{L^2} \le C_1 \, \frac{ R^k \|f_N\|^2 _{H^k _p}}{M^k} +
                    \frac{C_2}{N^k} \, \Big( \|f\| _{H^k _p} + \|Q^R (f_N,f_N)\|_{H^k _p} \Big).
  \end{multline*}
 \end{theorem}
\begin{proof}[Proof of Theorem~\ref{theo:consist}]
By triangular inequality
 \begin{multline*}
 \| Q^R (f,f) -  \mathcal{P}_N Q^{R,M} (f_N,f_N) \|_{L^2} \le
    \|  \mathcal{P}_N  \big( Q^R (f_Nf_N) - Q^{R,M}(f_N,f_N) \big) \|_{L^2} \\
    + \| Q^R (f,f) -  \mathcal{P}_N Q^{R} (f_N,f_N) \|_{L^2}.
 \end{multline*}
The first term on the right-hand side is controlled by
Lemma~\ref{lem:rect}
 \begin{multline*}
 \|  \mathcal{P}_N  \big( Q^R (f_N,f_N) - Q^{R,M}(f_N,f_N) \big) \|_{L^2} \le
        \| Q^R (f_N,f_N) - Q^{R,M}(f_N,f_N) \|_{L^2} \\
     \le C_1 \, \frac{ R^k \|f_N\|^2 _{H^k _p (\cal{D_\pi})}}{M^k}.
 \end{multline*}
The second term in the right-hand side is controlled by
Lemma~\ref{lem:PaRu}, which concludes the proof.
\end{proof}

Now let us focus briefly on the macroscopic quantities. In fact
here no additional error (related to $M$) occurs, compared with
the usual spectral method, since the approximation of the
collision operator that we are using is still conservative.
First with Lemma~\ref{lem:rect} at hand one can establish the estimate
 \[ \|Q^{R,M} (g,f) \|_{L^2} \le C \, \|g\|_{H^d _p} \|f\|_{H^d _p}, \]
for a constant uniform in $M$.
Then following the method of~\cite[Remark~5.4]{PaRu:spec:00}
and using this estimate we obtain the following spectral accuracy result
 \begin{multline*}
 \big| \langle Q^{R,M}(f,f), \varphi \rangle - \langle \mathcal{P}_N Q^{R,M}(f_N,f_N), \varphi \rangle \big|_{L^2}
 \le \frac{C_3}{N^k} \, \|\varphi\|_{L^2} \, \Big( \|f\| _{H^{k+d} _p} + \|Q^{R,M} (f_N,f_N)\|_{H^k _p} \Big)
 \end{multline*}
where $\varphi$ can be replaced by $v,|v|^2$. Indeed there is no
need to compare the momenta of $\mathcal{P}_N Q^{R,M}(f_N,f_N)$ with
those of $Q^R(f,f)$ since $Q^{R,M}$ is also conservative,
and so they can be compared directly to those of $Q^{R,M}$. Thus the error on
momentum and energy is independent on $M$ and is spectrally small according to $N$
even for very small value of the parameter $M$.


\subsection{Implementation of the algorithm}

The final spectral scheme depends on the three parameters $N$,
$R$, and $M$. The only conditions on these parameters is the
no-aliasing condition that relates $R$ and the size of the box $T$
(here $\pi$). A detailed study of the influence of the choices of
$N$ and $R$ has been done in \cite{PaRu:spec:00}. Here we are
interested only in the influence of $M$ over the computations,
since $M$ controls the computations speed-up.

The method of the previous subsections yields a decomposition of the
collision operator, which after
projection on $\ens{P}^N$ gives the following decomposition
 \begin{equation}\label{eq:opdecproj}
 \mathcal{P}_N Q^{R,M} = \sum_{p=1} ^{M^{d-1}} \mathcal{P}_N Q^{R,M} _p.
 \end{equation}
Each $\mathcal{P}_N Q^{R,M} _p$ can be computed with a cost $O(N^d
\log_2 N)$. Thus for a general choice of $M$ and $N$ we obtain the
cost $O(M^{d-1} N^d \log_2 N)$. The
decomposition~\eqref{eq:opdecproj} is completely
parallelizable and thus the cost can be strongly reduced on a
parallel machine (theoritically up to $O(N^d \log_2 N)$). One just has
to make independent computations for the $M^{d-1}$ terms of the
decomposition.

Moreover the formula of decomposition is naturally adaptive
(that is the number $M$ can be made space dependent), which can be
quite useful in the inhomogeneous setting, where some regions
deserve less accuracy than others. Since it relies on the
rectangular formula, whose adaptivity property is well known, one
can easily double the number of directions $M$ if needed, without
computing again those points already computed.

Finally the decomposition can be
also interesting from the storage viewpoint, as the classical spectral method
requires the storage of a $N^d \times N^d$ matrix
whereas our method requires the storage of $2 M^{d-1}$ vectors of size $N^d$.
In dimension $2$ the classical method requires a storage of order $O(N^4)$ and our
method requires a storage of order $O(M N^2)$. In dimension $3$
the classical method requires a storage of order $O(N^4)$ (thanks to the
symmetries of the matrix of kernel modes, see~\cite{FiRu:pre}), and our
method requires a storage of order $O(M^2 \, N^3)$.

As a numerical example we report the results obtained in the case
of space homogeneous two-dimensional Maxwellian molecules using as a
comparison the exact analytic solution (see \cite{PaRu:spec:00}).
The results for the relative $L_1$ norm of the error at time
$t=0.01$ are reported in Table \ref{table1}.

\begin{table}
\begin{center}
\begin{tabular}{|c|c|c|c|c|}
\hline
$N$& M=2& M=4& M=8& M=16\\
\hline  32 & 2.129E-4 & 1.993E-05 & 2.153E-05  & 2.262E-5 \\
\hline  64 & 2.109E-4 & 7.122E-10 & 6.830E-10  & 6.843E-10 \\
\hline 128 & 2.112E-4 & 3.116E-12 & 3.117E-12  & 3.117E-12 \\
\hline
\end{tabular}
\vskip .5cm \caption{Relative $L_1$ norm of the error for
different values of $N$ and $M$ for the fast spectral method. }
\label{table1}
\end{center}
\end{table}

Although further extensive testing is necessary, the results are
very promising and seem to indicate a very low influence of the
number of directions over the accuracy of the scheme. For $M=2$
the angle error dominates, but as soon as $M=4$ the error in $N$
is dominating. Note that the number of angle directions will
indirectly influence the aliasing effect trough the slight change in
the relaxation times. This may explain the slight error variations
that we observe taking $M\geq 4$.

Finally, in view of space non homogeneous computations, we will
have the additional advantage of taking a larger number of
gridpoints without increasing too much the computational cost,
thus allowing the computations of flows at larger Mach number
compared to conventional deterministic schemes. Further numerical
results are under development and will be presented in the
work \cite{FiMoPa:03}.

\section{Conclusions}
We have presented a deterministic way for computing the Boltzmann
collision operator with fast algorithms, for a class of interactions
which includes the case of hard spheres in dimension $3$.
The method is based on a
Carleman-like representation of the operator that allows to
express it as a combination of convolutions (this is trivially
true for the loss part but it is not trivial for the gain part). A
suitable periodized truncation of the operator is then used to
derive new spectral methods 
computable with a high speed up in computation times.
This brings the overall cost in dimension $d$ to
$O(M^{d-1}N^d\log_2 N)$ where $N$ is the number of velocity
parameters and $M$ the number of angular directions in each
angular coordinate.
Consistency and accuracy of the proposed
schemes are also presented, and it is shown to be spectrally accurate.
Moreover the error on the momentum
and energy is spectrally small and independent of the
value of the speed-up parameter $M$.
First numerical results seem to
indicate the validity and the flexibility of the present approach
that, to our opinion, will make deterministic schemes much more
competitive with Monte Carlo methods in several situations.

\appendix
\section*{Appendix: Remarks on admissible collision kernels and an extension to the ``non-decoupled'' case}

Let us study the cases where the assumption~\eqref{eq:decoup} is
satisfied. For hard spheres in dimension $3$, or Maxwellian
molecules in dimension $2$, one has the equation~\eqref{eq:decoup}
with $a=b=1$. Formally for the Coulomb potential in dimension $3$,
we have
 \[ B(\theta,|u|) = |u|^{-3} \, \sin^{-4} (\theta/2), \]
and thus, thanks to formula~\eqref{eq:Btilde}
 \[ \tilde{B}(x,y) = 2^{d-1} \, |x|^{-4}. \]

This suggests, in dimension $3$, to consider the following family of ``variable
hard sphere'' collision kernels
 \begin{equation}\label{eq:VHSct}
 B_\gamma(\theta,|u|) = \sin^{\gamma-1} (\theta/2) \,
 |u|^{\gamma}.
 \end{equation}
Indeed simple computations give
 \[ \tilde{B_\gamma}(x,y) = 2^{d-1} \, |x|^{\gamma-1} \]
and thus they satisfy the decoupling assumption~\eqref{eq:decoup}. 
In the case where $\gamma \in (-2,1]$ the angular part of the
collision kernel remains integrable. On the contrary, for $\theta
\sim 0$, the equivalent derived from the physical non explicit
formula in~\cite{Cer:88} for inverse-power laws kernels (for a
potential $1/|d|^{n-1}$ with $n$ such that $\gamma = (n-5)/(n-1)$) is
of the form
 \begin{equation}\label{eq:VHS} \nonumber
 B^{\mbox{\scriptsize{exact}}} _\gamma (\theta,|u|) \sim_{\theta \sim 0} \sin^{\frac{\gamma-5}2} (\theta/2) \, |u|^{\gamma}
 \end{equation}
with $\gamma \in [-3,1)$. It is therefore always non-integrable for $\theta \sim 0$.


The model~\eqref{eq:VHSct} coincides with the hard spheres model for $\gamma
=1$ and, formally, coincides with the kernel of the Coulomb potential for $\gamma =
-3$. Moreover for $\gamma \in (-2,1]$ ({\em i.e.} hard potentials and the so-called
moderately soft potentials) it remains integrable for $\theta \sim 0$.
Thus it seems quite reasonable to consider it as a model for cutoff hard
and moderately soft potentials, as well as hard spheres.

In dimension $2$ the same arguments and computations lead to the following cutoff hard and
moderately soft potentials model
 \begin{equation}\label{eq:VHSct2} \nonumber
 B_\gamma (\theta,|u|) = \sin^\gamma (\theta/2) \, |u|^{\gamma}
 \end{equation}
valid for $\gamma \in (-3,1]$, which coincides with the case of Maxwellian molecules for $\gamma =0$.

For the spectral method the other situation where one obtains naturally a fast algorithm is
the case where collisions concentrate on the grazing part: see~\cite{PaRuTo:fastFPL:00}
and~\cite{PaToVi:03} for a fast algorithm to compute the Fokker-Planck-Landau collision operator, which is the
limit of the Boltzmann collision operator in the grazing collision limit. In this case indeed one
of the two variables $x$ or $y$ of the representation~\eqref{eq:repR}
disappears in the limit process, which ``decouples'' the kernel modes.
Thus it may be possible to construct fast algorithms for non-cutoff models by splitting the collision
operator into a cutoff part treated by the method presented in this paper, and a non-cutoff part restricted
to very small deviation angles, which would be close to the grazing collision limit and thus could
be computed by the fast algorithm of~\cite{PaRuTo:fastFPL:00,PaToVi:03}.


\bigskip

\noindent {\bf{Acknowledgments.}}
Both authors thank Francis Filbet for the numerical results of
Table 1. Support by the European network HYKE, funded by the EC as
contract HPRN-CT-2002-00282, is acknowledged.

\medskip

\bibliographystyle{acm}

\bibliography{bibMPfast}

\begin{flushleft} \signcm \end{flushleft}
\vspace*{-45mm}
\begin{flushright} \signlp \end{flushright}

\end{document}